\documentclass[12pt,reqno,oneside]{amsart}
\usepackage{amsmath,amsthm,amsfonts,amssymb}
\usepackage{indentfirst}
\usepackage{url}
\usepackage{hyperref}
\usepackage{graphicx}

\usepackage{caption}
\usepackage{subcaption}

\theoremstyle{plain}
\newtheorem{theorem}{Theorem}[section]
\newtheorem{cor}[theorem]{Corollary}

\theoremstyle{definition}

\newtheorem{obs}[theorem]{Remark}

\numberwithin{equation}{section}

\begin{document}

\baselineskip=18pt

\title[The fitness of the strongest individual]{The fitness of the strongest individual in the subcritical GMS model}

\author[Carolina Grejo]{Carolina Grejo}
\address[F. Machado, A. Rold\'an, C. Grejo]{Statistics Department, Institute of Mathematics and Statistics, University of S\~ao Paulo, CEP 05508-090, S\~ao Paulo, SP, Brazil.}
\email{carolina@ime.usp.br}
\author[F\'abio Machado]{F\'abio Machado}
\email{fmachado@ime.usp.br}
\thanks{Carolina Grejo was supported by CNPq (141965/2014-2), F\'abio Machado by CNPq (310829/2014-3) and Fapesp(09/52379-8) and Alejandro Roldan by CNPq (141046/2013-9).}
\author[Alejandro Rold\'an-Correa]{Alejandro Rold\'an-Correa}
\email{roldan.alejo@gmail.com}

\date{\today}

\begin{abstract} 
{We derive the strongest individual fitness distribution on a variation for a species 
survival model proposed by Guiol, Machado and Schinazi~\cite{GMS11}. We point out
to the fact that this distribution relies on the Gauss hypergeometric function and when
$p=\frac{1}{2}$ on the hypergeometric function type I distribution.} 
\end{abstract}

\maketitle

\section{Introduction}
\label{S: Introduction}

We consider a discrete time model  beginning from an empty set. At each time $n\geq 1$, a new species is born with probability $p$ or there is a death (if the system is not empty) with probability $q=1-p$. Let $X_n$ be the total number of species at time $n$.  $X_n$ is a random walk on ${\mathbb{Z}}_+$ that jumps to right with probability $p$ and jumps to left with probability $q$. When $X_n$ is at 0 the process jumps to 1 with probability $p$ or stays at 0 with probability $1-p$. We assign a random number to each new species. This number has a uniform distribution on $[0,1]$. We think of this number as a fitness associated to each species. These random numbers are independent to each other. When a death occurs, the individual with lowest fitness dies. This model, latter denominated GMS model, was first proposed and studied in Guiol {\it et al}~\cite{GMS11}. Some interesting variations were further studied in Guiol {\it et al}~\cite{GMS13}, Ben Ari {\it et al}~\cite{BMR11} and Skevi and Volkov~\cite{SV12}.

In Guiol {\it et al} \cite{GMS11} it is shown that there is a sharp phase transition for $p>1/2$. 
For $R_n$, the set of species with fitness higher than $f_c=\frac{1-p}{p}$ at time $n$ 
approachs an uniform distribution in the following sense. For $f_c < a < b <1$
\[ \lim_{n \to \infty} \frac{|R_n \cap (a,b)|}{n} = p(b-a)  \ \text{ a.s.}\]
On the other hand every specie born with fitness less than $f_c$  disappear after a finite (random) time. 
The set of species present in the system whose fitness is smaller than $f_c$ becomes empty infinitely many times.

Here we focus on the case $p\leq1/2$ in order to understand 
better the dynamics of this model. In this case, the process $X_n$ is recurrent and the system becomes empty infinitely many times. Therefore it is not interesting to study the distribution of the fitness of the species which are alive on the system in the long run. An interesting point is to study the distribution of the fitness of the strongest individual on each excursion between the epochs when the system becomes empty. 

We propose a variation for the GMS model by considering that each time the system becomes empty, a set of $m$ individuals are introduced with independent set of fitness. This variation is meant to reinforce competition among species before the system becomes empty again.

\section{Results}

We deduce explicitly the distribution of the fitness of the strong\-est individual on excursions between the epochs when the system becomes empty. The last individual to die before the system becomes empty is the strongest on that excursion because the first ones to die are those individuals with the smallest fitness.

Observe that some excursions may have length 2. When this happens, the individual who is born, dies right away without competing with any other individual. To ensure that each excursion has competition among individuals in a sort of natural selection process, we introduce a change-over on the model: Each time after the system becomes empty, $m$ independent new species are placed on the system (instead of just 1) with probability $p$, or the system stays empty
with probability $1 - p$. We denote this variation by GMS($m$). In this set up GMS(1) is the original model.

\begin{figure}
\centering
\begin{subfigure}{.5\textwidth}
  \centering
  \includegraphics[scale=0.3]{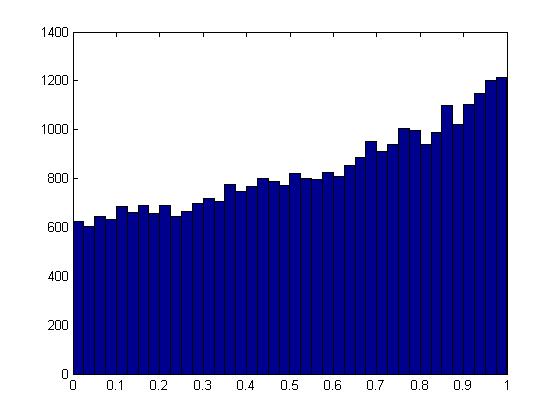}
  \caption{GMS(1)}
  \label{fig:sub1}
\end{subfigure}%
\begin{subfigure}{.5\textwidth}
  \centering
  \includegraphics[scale=0.3]{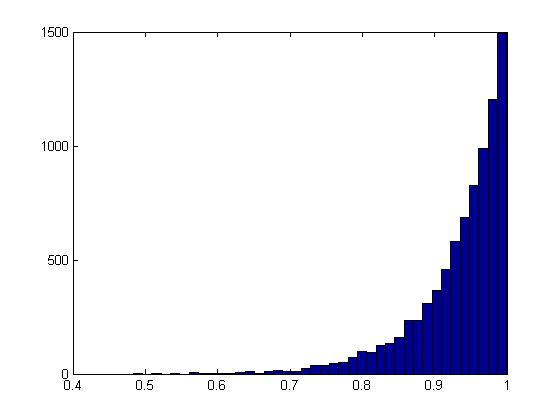}
  \caption{GMS(10)}
  \label{fig:sub2}
\end{subfigure}
\caption{Histograms of the fitnesses of the strongest individual on GMS(m) after 200,000 births and deaths for $p=1/4$.}
\label{fig:test}
\end{figure}

Figures \ref{fig:sub1} and \ref{fig:sub2}
show the role of the competition on the distribution
of the fitness of the strongest individual on each excursion. Short excursions are more commom on GMS(1) than on 
GMS(10). That behaviour favors individuals with lower fitnesses to be the strongest ones. Competition 
introduced in GMS(10) avoids that.

The next result computes the fitness distribution of the strongest individual to die right before the system becomes 
empty on GMS($m$) model. It is shown in terms of the hypergeometric function of Gauss (see Luke~\cite{Luke}). This function is denoted by ${_2F_1}(a,b;c;z)$, namely,
\begin{equation}\label{HG}
{_2F_1}(a,b;c;z)=\sum_{k\geq 0}\frac{(a)_k (b)_k}{(c)_k}\frac{z^k}{k!}, \quad |z|<1,
\end{equation}

\noindent where $a$, $b$, $c$, are real numbers with $c\neq 0,-1,-2,\dots$, and  $(a)_k$ is the coefficient Pochhammer, namely, 
$$(a)_k=a(a+1)\cdots(a+k-1) \quad (a)_0=1.$$

\begin{theorem}\label{T:GMSm} Let  $p\leq 1/2$ and  $Z_m$ be the fitness of the strongest individual before the system 
becomes empty on  GMS($m$) model. Then  $Z_m$ is a random variable with distribution 
$$\mathbb{P}[Z_m\leq t]= (qt)^m {_2F_1}\left(\frac{m}{2},\frac{m+1}{2};m+1;4pqt\right) \hskip1cm 0\leq t <1.$$
\end{theorem}

\begin{cor}\label{T:GMS}  Let  $p\leq 1/2$ and $Z$ be the fitness of the strongest 
individual before the system becomes empty on  GMS($1$) model. Then
$$\mathbb{P}[Z\leq t]=\frac{1-\sqrt{1-4pqt}}{2 p}, \hskip1cm 0\leq t <1.$$
For $p= 1/2$, $Z$ follows a Beta distribution $B(1,1/2).$
\end{cor}


\begin{obs} By Theorem \ref{T:GMSm}  we have $Z_m$ density probability function is
\begin{eqnarray}\label{hipertipo1}
f_m(t) 
&=&  \frac{d}{dt}\left[(qt)^m {_2F_1}\left(\frac{m}{2},\frac{m+1}{2};m+1;4pqt\right)\right]\nonumber\\
&=&   mq^m t^{m-1} {_2F_1}\left(\frac{m}{2},\frac{m+1}{2};m;4pqt\right)      
\end{eqnarray}
where the last line have been obtained by using Abramowitz and Stegun \cite[Eq. 15.2.4]{Abramowitz}.  When  $p=q=1/2$, the distribution of $1-Z_m$ is known as \textit{hypergeometric function type I distribution} (see Gupta and Nagar \cite[p. 298]{Nagar}). 
\end{obs}

\begin{cor}\label{C:GMSm}  $\mathbb{E}[Z_m]=1-\displaystyle\frac{q^m}{m+1} \ {_2F_1}\left(\frac{m}{2},\frac{m+1}{2};m+1;4pq\right)$
\end{cor}

\begin{obs} Considering Corollary~\ref{C:GMSm} when $p=1/2$, by using the following equality  (see Gradshteyn and Ryzhik \cite[Eq. 7.512.11]{Gradshteyn}) $${_2F_1}(\alpha,\beta;\gamma;1)=\frac{\Gamma(\gamma)\Gamma(\gamma-\alpha-\beta)}{\Gamma(\gamma-\alpha)\Gamma(\gamma-\beta)}$$ we have that

\begin{eqnarray}
\mathbb{E}[Z_m]&=&1-\frac{q^m}{m+1} {_2F_1}\left(\frac{m}{2},\frac{m+1}{2};m+2;1\right) \nonumber \\ 
&=&1-\frac{1}{2^m(m+1)}\left[\frac{\Gamma(m+2)\Gamma(3/2)}{\Gamma(\frac{m+3}{2}+\frac{1}{2})\Gamma(\frac{m+3}{2})}\right] \nonumber \\
&=&1-\frac{\sqrt{\pi} \ m!}{2^{m+1}\Gamma(\frac{m+3}{2}+\frac{1}{2})\Gamma(\frac{m+3}{2})} \nonumber
\end{eqnarray}

\noindent where the last line has been obtained by using $\Gamma(3/2)=\sqrt{\pi}/2$. Now, using the duplication formula, namely, 

\begin{equation} 
\Gamma(2z)=\frac{\Gamma(z+\frac{1}{2})\Gamma(z)}{2^{1-2z}\sqrt{\pi}} \nonumber
\end{equation}
\noindent we get 
\begin{equation} 
\mathbb{E}[Z_m]=1-\displaystyle\frac{2}{(m+1)(m+2)} \nonumber
\end{equation}
\end{obs}

\section{Proofs}

\begin{proof}[Proof Theorem \ref{T:GMSm}]
For $n=0,1, \dots$ we define
\[ \tau_n = \inf \{l \ge 1 : X_{n+l} = 0, X_n=0 \}. \]

In words $\tau_n$ is the length of a excursion from 0 to 0. As the process $X_n$ is homogeneous,
the distribution of $\tau_n$ does not depend on $n$ so we consider the random variable 
$\tau := \tau_0$. Besides, as $p\leq 1/2$ we have that $\mathbb{P}[\tau <\infty]=1$ and
$$\mathbb{P}[\tau =k+1]=\mathbb{P}[T_{-m}=k]=\frac{m}{k}{k \choose \frac{k-m}{2}} p^{(k-m)/2}q^{(k+m)/2}, \ k\geq m, \ k+m \text{ even}, $$
where $T_{-m}$  is the time of the first visit to  $-m$ for a random walk on $\mathbb{Z}$ beginning at 0. (See Bhattacharya and Waymire \cite{Batachyara})\\

\noindent If $\tau=k+1,$ we see along that excursion, extra $\frac{k-m}{2}$ births and $\frac{k+m}{2}$ deaths. 
The last death corresponds to the individual with the strongest fitness among all $\frac{k+m}{2}$ that were born. Hence,

$$\mathbb{P}[Z_m\leq t]=\sum_{k=m}^\infty \mathbb{P}[\tau=k+1]\mathbb{P}[\max(Y_1,...,Y_{\frac{m+k}{2}})\leq t],$$

\noindent 
where  $Y_1,...,Y_{\frac{m+k}{2}}$ are i.i.d. uniform random variables on  $[0,1].$ Therefore,

$$\begin{array}{ccl}
\mathbb{P}[Z_m\leq t]&=&\displaystyle\sum_{k=m}^\infty \frac{m}{k}{k \choose \frac{k-m}{2}} p^{(k-m)/2}q^{(k+m)/2}t^{(k+m)/2} \textbf{1}_{\{m+k \text{ even}\}} \\ \\
&=&\displaystyle\sum_{l=m}^\infty \frac{m}{2l-m}{2l-m \choose l} p^{l-m}q^{l}t^{l} \hskip0.5cm (k+m=2l, l\geq m)\\  \\
&=&\displaystyle\sum_{j=0}^\infty \frac{m}{m+2j}{m+2j \choose m+j} p^{j}q^{m+j}t^{m+j} \hskip0.5cm (l=m+j)\\ \\
&=&(qt)^m \displaystyle\sum_{j=0}^\infty \frac{m}{m+2j}{m+2j \choose m+j} (pqt)^{j}\\
&=&(qt)^m \displaystyle\sum_{j=0}^\infty \frac{(m)_{2k}}{(m+1)_k} \frac{(pqt)^{j}}{k!}\\ \\
&=&(qt)^m {_2F_1}\left(\frac{m}{2},\frac{m+1}{2};m+1;4pqt\right) \\ \\
\end{array}$$
where the last line has been obtained by using $(a)_{2k}=\left(\frac{a}{2}\right)_k\left(\frac{a+1}{2}\right)_k  2^{2k}$  and the definition of Gauss hypergeometric function.
\end{proof}

\begin{proof}[Proof Corollary \ref{T:GMS}]
It is a particular case of Theorem \ref{T:GMSm}  when $m=1$. In this situation 
\begin{eqnarray}
qt \ {_2F_1}\left(\frac{1}{2},1;2;4pqt\right)
&=&  qt \sum_{k\geq 0}\frac{(1/2)_k (1)_k}{(2)_k}\frac{(4pqt)^k}{k!} \nonumber\\ 
&=&  \frac{1}{p} \sum_{k\geq 0} \frac{(2k)!}{(k+1)! \ k!}(pqt)^{k+1} \nonumber\\
&=&  \frac{1-\sqrt{1-4pqt}}{2p}\nonumber
\end{eqnarray}
where the last line has been obtained by using $(1)_k(1/2)_k= 2^{-2k} (2k)!$ and the result given in Prudnikov {\it et al} \cite[Eq. 5.2.13.8]{Prud}.
\end{proof}

\begin{proof}[Proof Corollary \ref{C:GMSm}]

$$\begin{array}{ccl}
\mathbb{E}[Z_m]&=& \int_0^1\mathbb{P}[Z_m>t] \ dt \\ \\
&=&1-q^m\int_0^1 t^m {_2F_1}\left(\frac{m}{2},\frac{m+1}{2};m+1;4pqt\right) \ dt \\ \\
&=&1-\frac{q^m}{m+1} {_2F_1}\left(\frac{m}{2},\frac{m+1}{2};m+2;4pq\right)\\ \\
\end{array}$$
where the last line has been obtained by using the result given in Gradshteyn and Ryzhik \cite[Eq. 7.512.11]{Gradshteyn}.
\end{proof}

\noindent
\textbf{Acknowledgements:} The authors are thankful to Daniel Valesin and Rinaldo Schinazi for helpful discussions about the model. Thanks are also due to the anonymous referee for his/her
constructive comments, leading to an improved presentation.

\end{document}